\documentclass{amsart}

\newtheorem{tm}{Theorem}
\newtheorem{df}{Definition}
\newtheorem{lm}{Lemma}
\newtheorem{kor}{Corollary}

\newtheorem{conj}{Conjecture}

\theoremstyle{remark}
\newtheorem{remark}{Remark}

\usepackage[usenames,dvipsnames,svgnames,table]{xcolor}
\usepackage{comment}

\begin{document}

\title[The extensibility of the Diophantine triple $\{2,b,c\}$]{THE EXTENSIBILITY OF THE DIOPHANTINE TRIPLE $\{2,b,c\}$ }

\author{NIKOLA AD\v{Z}AGA, ALAN FILIPIN AND ANA JURASI\'{C}}

\begin{abstract}The aim of this paper is to consider the extensibility of the Diophantine triple $\{2,b,c\}$, where $2<b<c$, and to prove that such a set cannot be extended to an irregular Diophantine quadruple. We succeed in that for some families of $c$'s (depending on $b$). As corollary, for example, we prove that for $b/2-1$ prime, all Diophantine quadruples $\{2,b,c,d\}$ with $2<b<c<d$ are regular. 

\end{abstract}

\maketitle


\section{INTRODUCTION}

A set consisting of $m$ positive integers such that the product of any two
of its distinct elements increased by 1 is a perfect square is called a Diophantine $m$-tuple. There is long history of finding such sets. One of the questions of interest, which various mathematicians try to solve, is how large those sets can be. Very recently He, Togb\'{e} and
Ziegler \cite{HTZ} proved the folklore conjecture that there cannot be 5 elements in Diophantine $m$-tuple, i.e. $m<5$. However, there is also a stronger version
of that conjecture that is still open, which states that every Diophantine triple can
be extended to a quadruple with a larger element in a unique way
(see \cite{duje1}):

\begin{conj}\label{conj:1}
If $\{a,b, c, d\}$ is a Diophantine quadruple of integers and $d>\rm{max}\{\textit{a,b,c}\}$,
then $d = d_{+}=a+b+c+2(abc+\sqrt{(ab+1)(ac+1)(bc+1)}).$
\end{conj}

There are a lot of results supporting this conjecture. The history of the problem, with recent results and up-to-date references can be found on the webpage \cite{web}.

We have the following definitions (see \cite{gib-dodano}). Let
$\{a,b,c\}$ be a Diophantine triple such that
\begin{eqnarray}\label{prva}
ab+1=r^{2},\
ac+1=s^{2},\ bc+1=t^{2},
\end{eqnarray}
where $r,s,t\in\mathbb{N}$.
\begin{df}\label{def1}
A Diophantine triple $\{a,b,c\}$ is called regular if
\begin{equation} \label{zv1}
(c-b-a)^{2}=4(ab+1).
\end{equation}
\end{df}Equation \eqref{zv1} is obviously symmetric under permutations of $a$, $b$,
$c$. From \eqref{zv1}, by \eqref{prva}, we get
\begin{eqnarray*}c=c^{\pm}=a+b \pm 2r, \end{eqnarray*}
\begin{eqnarray*}ac^{\pm} + 1 = (a\pm r)^2, \, bc^{\pm} + 1 = (b\pm r)^2. \end{eqnarray*} \begin{df}\label{def2}
A Diophantine quadruple $\{a,b,c,d\}$ is called regular if
\begin{eqnarray}\label{drugo}
(d+c-a-b)^{2}=4(ab+1)(cd+1)
\end{eqnarray} or, equivalently, if
\begin{equation}\label{d-plus}d=d_{\pm}=a+b+c+2(abc\pm rst).\end{equation}
\end{df}By (\ref{prva}) and (\ref{d-plus}), we have \begin{equation}\label{d-st-reg}ad_{\pm} + 1 = (at\pm rs)^2, \, bd_{\pm} + 1 = (bs\pm rt)^2, \, cd_{\pm} + 1 = (cr\pm st)^2. \end{equation}Equation \eqref{drugo} is symmetric under permutations of $a$,
$b$, $c$, $d$. An irregular Diophantine quadruple is one that is not
regular. It is known from \cite{AHS} that every Diophantine pair $\{a,b\}$ can be extended to a regular Diophantine quadruple: \begin{eqnarray*}\{a,b,a+b\pm 2r,4r(a\pm r)(b\pm r)\}.\end{eqnarray*}

During the second conference on Diophantine $m$-tuples and related problems, that took place at Purdue University Northwest, we mentioned the following result \cite[Theorem 4]{ADKT}. If $\{2,b,c,d\}$ is a regular Diophantine quadruple, then the Diophantine triple $\{b,c,d\}$ is also a $D(n)$-set for two distinct $n$'s with $n\neq 1$ (which means that $bc+n$, $bd+n$ and $cd+n$ are perfect squares). We have realized that this result would be even more elegant if we could drop the word ``regular''. On the other hand, in \cite{HPTY}, Conjecture \ref{conj:1} was proven when $a=1$, so it makes sense to see if the method they used can prove this conjecture for different values of $a$ before attempting to generalize it, because that would probably be very difficult.

For $a=2$, by (\ref{prva}), we have \begin{equation} \label{k}b=2k(k+1),\ \ r=2k+1, \end{equation}for $k\in\mathbb{N}$. Hence, $b=4(1+...+k)$, where $k\in\mathbb{N}$. We can notice that $b$ is always even. We take that $b>4000$ because otherwise the triple $\{2,b,c\}$, with $2<b<c$, cannot be extended to an irregular quadruple by \cite[Lemma 3.4]{CFF}. We will use the condition $b>4000$ through the paper. Since $1^2+2^2+...+k^2 =\frac{1}{6}(2k+1)k(k+1)$, it is interesting to observe that $\frac{br}{12}=P_k$, where $P_k$ is a square pyramidal number for $k\in\mathbb{N}$. The problem of extending a Diophantine pair $\{2,b\}$ to a Diophantine triple $\{2,b,c\}$ can be reduced to solving the Pellian equation \begin{equation}\label{Pell}2t^2-bs^2=2-b,\end{equation}and then taking $c=\frac{s^2-1}{2}$. We will describe the set of solutions of the equation (\ref{Pell}), by following the arguments of Nagell \cite{Nagel} and Dujella \cite{duje-Pell}. Equation (\ref{Pell}) is equivalent to the Pellian equation \begin{equation}\label{Pell-D}t^2-\frac{b}{2}s^2=1-\frac{b}{2}.\end{equation}Such an equation has infinitely many solutions divided into classes\footnote{Two solutions $t+s\sqrt{b/2}$ and $t'+s'\sqrt{b/2}$ of the equation (\ref{Pell-D}) belong to the same class if and only if $tt'\equiv (b/2)ss'\ ({\rm mod}\ (1-b/2))$ and $ts'\equiv t's\ ({\rm mod}\ (1-b/2))$.} of solutions. Among all elements of one class, we choose $t_0+s_0\sqrt{b/2}$, where $s_0$ has the lowest nonnegative value among all elements  $t+s\sqrt{b/2}$ of the same class. Such a solution is called a fundamental solution of the equation (\ref{Pell-D}). Notice also that $|t_0|$ has the lowest possible value in the class. By the arguments of \cite[Theorem 108 a]{Nagel}, equation (\ref{Pell-D}) has finitely many fundamental solutions. Hence, it has finitely many classes of solutions. There are at most $2^{\omega(1-b/2)}$ classes of solutions with $\rm{gcd}(\textit{t},\textit{s})=1$, where $\omega(1-b/2)$ denotes the number of distinct prime factors of $1-b/2$. Although we will focus on the situation when there is only one class of solutions, with fundamental solutions $(t_0,s_0)=(\pm 1,1)$, let us mention here that this is not always the case. If $k$ is of the form $k=g^2-2$, then \eqref{Pell-D} becomes $t^2-(g^4-3g^2+2)s^2=-g^4+3g^2-1$ and this equation has fundamental solution $(g^3-g^2-2g+1, g-1)$, which differs from $(\pm 1, 1)$ when $g\neq 2$.

All elements of one class of solutions of the equation (\ref{Pell-D}) can be obtained from a fundamental solution by multiplication with a power of the minimal solution in positive integers for the associated Pellian equation $t^2-\frac{b}{2}s^2=1$. Therefore, all positive solutions $(t,s)$  to equation (\ref{Pell}) which belong to the same class are given with (see \cite[Lemma 1]{duje-Pell}) \begin{eqnarray}\label{Pell-general}t\sqrt{2}+s\sqrt{b}=(t_0\sqrt{2}+s_0\sqrt{b})(r+\sqrt{2b})^\nu,\end{eqnarray}where $(t_0,s_0)$ is a fundamental solution to the equation (\ref{Pell}) and $\nu\geq 0$. Since we have fundamental solutions $(t_0,s_0)=(\pm 1,1)$, then if those are only fundamental solutions (for example if $b/2-1$ is a prime, but this is not the only case), all positive solutions to the equation (\ref{Pell}) are given with $$(t,s)=(t_{\nu}^{\pm},s_{\nu}^{\pm}),$$ obtained by the recurrent relations for two different signs $\pm 1$:
\begin{eqnarray}\label{rekurzija1}
t_{0}&=&\pm 1,\ \ t_{1}=b\pm r,\ \ t_{\nu +2}=2rt_{\nu +1}-t_{\nu},\\
s_{0}&=&1,\ \ s_{1}=r\pm 2,\ \
s_{\nu +2}=2rs_{\nu +1}-s_{\nu},\label{rekurzija2}\end{eqnarray}for $\nu\geq 0$. Since for $\nu = 0$ we obtain $c = 0$, we have to consider the sequence $c=c_{\nu}^\pm$, with $2c_\nu^\pm+1=(s_\nu^\pm)^2$, $bc_\nu^\pm+1=(t_\nu^\pm)^2$, for $\nu\geq 1$. By \cite[Theorem 1.4. (4)]{c-f-m}, we need to consider only $1\leq \nu \leq 3$ i.e.
\begin{eqnarray*}c_1^\pm&=&2+b\pm 2r,\\
c_2^\pm&=&8b(2+b\pm 2r)+4(2+b\pm r)=4r(r\pm 2)(b\pm r),\\
c_3^\pm&=&64b^2(2+b\pm 2r)+16b(6+3b\pm 4r)+3(6+3b\pm 2r).\end{eqnarray*} We observe the extensibility of the triple $\{2,b,c_\nu^\pm\}$, where $1\leq \nu \leq 3$. Since $bc_1^-+1=(b-r)^2<b^2$, it follows that $c_1^-<b$. In all other cases $b<c_\nu^\pm$. Our main result in this paper is the following theorem:
\begin{tm}\label{Tm1}
The triple $\{2,b,c_\nu^\pm\}$, for $\nu\in\mathbb{N}$, cannot be extended to an irregular Diophantine quadruple $\{2,b,c_\nu^\pm,d\}$, where $d>c_\nu^\pm$.\end{tm}

Theorem \ref{Tm1} allows us to derive the following statement from the previous observations.

\begin{kor}If $\frac{b}{2}-1$ is a prime, then every Diophantine quadruple $\{2,b,c,d\}$, with $2<b<c<d$,  is regular.\label{Cor1}\end{kor}

In order to prove Theorem \ref{Tm1}, we use methods described in \cite{HPTY} by He, Pint$\rm{\acute{e}}$r, Togb$\rm{\acute{e}}$ and Yang. In Section 2 we transform the problem of
extending a Diophantine triple $\{2,b,c\}$ to a
Diophantine quadruple $\{2,b,c,d\}$ into solving a system of
simultaneous Pellian equations, which furthermore transforms to finding
intersections of binary recurrent sequences. In the next two sections we finish our proofs using the standard methods, i.e. linear forms in three, respectively two, logarithms of algebraic numbers and Baker-Davenport reduction method.

\section{The system of simultaneous Pellian equations}

When trying to extend of a Diophantine triple $\{2,b,c\}$ to a quadruple $\{2,b,c,d\}$ we have to find $d,x,y,z\in\mathbb{N}$ such that\begin{eqnarray}\label{treca}2d+1=x^{2},\ \  bd+1=y^{2},\ \ cd+1=z^{2}.\end{eqnarray}Elimination of $d$ from (\ref{treca}) leads to the system of simultaneous Pellian equations\begin{eqnarray}\label{cetvrta}2z^{2}-cx^{2}&=&2-c,\\
bz^{2}-cy^{2}&=&b-c,\label{peta}\\
2y^{2}-bx^{2}&=&2-b.\label{peta'}\end{eqnarray}Each of equations (\ref{cetvrta}), (\ref{peta}) and (\ref{peta'}) has finitely many fundamental solutions $(z_0,x_0)$, $(z_1,y_1)$ and $(y_2,x_2)$, respectively. From these solutions, all solutions $(z,x)$, $(z,y)$ and $(y,x)$ of (\ref{cetvrta}),  (\ref{peta}) and  (\ref{peta'}), respectively, are, by \cite[Lemma 1]{duje-Pell}, given with
\begin{eqnarray}\label{deseto}
 z\sqrt{2}+x\sqrt{c}&=&(z_{0}\sqrt{2}+x_{0}\sqrt{c})(s+\sqrt{2c})^{m},\\
 \label{jedanaesto}
 z\sqrt{b}+y\sqrt{c}&=&(z_{1}\sqrt{b}+y_{1}\sqrt{c})(t+\sqrt{bc})^{n},\\
\label{jedanaesto'}
 y\sqrt{2}+x\sqrt{b}&=&(y_{2}\sqrt{2}+x_{2}\sqrt{b})(r+\sqrt{2b})^{l},
 \end{eqnarray}for integers $m,n,l\geq 0$. In any solution $(x,y,z)$ of the system (\ref{cetvrta}) - (\ref{peta'}), we have $z=v_{m}=w_{n}$, for some non-negative integers $m$ and $n$, where the sequences
$(v_{m})_{m\geq 0}$ and $(w_{n})_{n\geq 0}$ are obtained using (\ref{deseto}) and (\ref{jedanaesto}) and given by
\begin{eqnarray*}
v_{0}&=&z_{0},\ \ v_{1}=sz_{0}+cx_{0},\ \ v_{m+2}=2sv_{m+1}-v_{m},\\
w_{0}&=&z_{1},\ \ w_{1}=tz_{1}+cy_{1},\ \ 
w_{n+2}=2tw_{n+1}-w_{n}.\end{eqnarray*}Hence, we are solving the equation \begin{equation} \label{jednadzbica}v_m=w_n\end{equation}in $m,n \geq 0$. By \cite[Theorem 2.1., Lemma 2.3.]{c-f-m}, it is enough to observe the cases:
\begin{enumerate}
\item [\textbullet]   if $v_{2m}=w_{2n}$, then $z_{0}=z_{1}= \pm1$ and $x_0=y_1=1$,
\item [\textbullet] if $v_{2m+1}=w_{2n+1}$, then $z_0=\pm t$, $z_1=\pm s$, $x_0=y_1=r$ and $z_0z_1>0$.
\end{enumerate}

In any solution $(x,y,z)$ of the system (\ref{cetvrta}) - (\ref{peta'}), we also have $y=W_{n}=q_{l}$, for some non-negative integers $n$ and $l$, where the sequences
$(W_{n})_{n\geq 0}$ and $(q_{l})_{l\geq 0}$ are obtained using (\ref{jedanaesto}) and (\ref{jedanaesto'}) and given by
\begin{eqnarray}\label{osma'}
W_{0}&=&y_{1},\ \ W_{1}=ty_{1}+bz_{1},\ \ W_{n+2}=2tW_{n+1}-W_{n},\\
q_{0}&=&y_{2},\ \ q_{1}=ry_{2}+bx_{2},\ \ \label{deveta'}
q_{l+2}=2rq_{l+1}-q_{l}.\end{eqnarray}We finally have $x=V_{m}=p_{l}$, for some non-negative integers $m$ and $l$, where the sequences
$(V_{m})_{m\geq 0}$ and $(p_{l})_{l\geq 0}$ are obtained using (\ref{deseto}) and (\ref{jedanaesto'}) and given by
\begin{eqnarray}\label{osma}
V_{0}&=&x_{0},\ \ V_{1}=sx_{0}+2z_{0},\ \ V_{m+2}=2sV_{m+1}-V_{m},\\
p_{0}&=&x_{2},\ \ p_{1}=rx_{2}+2y_{2},\ \ \label{deveta}
p_{l+2}=2rp_{l+1}-p_{l}.\end{eqnarray}From (\ref{osma'}) and (\ref{deveta'}), by induction, we get
\begin{eqnarray*}
W_{2n}&\equiv& y_{1}\ ({\rm mod}\ {b}),  \ \ W_{2n+1}\equiv ty_{1} ({\rm mod}\ {b}), \\
q_{2l}&\equiv& y_{2}\ ({\rm mod}\ {b}),  \ \ q_{2l+1}\equiv ry_{2}\ ({\rm mod}\ {b}).
\end{eqnarray*}By \cite[Lemma 1]{duje1}, $$1\leq |y_2|\leq\sqrt{\frac{(r-1)(b-2)}{4}}<\frac{b}{2}.$$ Therefore, there hold:
\begin{enumerate}
\item [\textbullet]  If $W_{2n}=q_{2l}$, then $1\equiv  y_{2}\ ({\rm mod}\ {b})$ so $y_2=1$.
\item [\textbullet]  If $W_{2n}=q_{2l+1}$, then $1\equiv  ry_{2}\ ({\rm mod}\ {b})$ so  $r\equiv  r^2y_{2}\ ({\rm mod}\ {b})$. By (\ref{prva}), $y_2\equiv  r\ ({\rm mod}\ {b})$. Since, for $b>4000$, $|y_2|<\frac{b}{2}$ and $r<\frac{b}{2}$, we get $y_2=r$. In this case, from (\ref{peta'}), we get $x_2^2=5$, which is not possible.
\item [\textbullet]  If $W_{2n+1}=q_{2l}$, then $rt\equiv  y_{2}\ ({\rm mod}\ {b})$. From (\ref{rekurzija1}), we conclude that $t=t_\nu^\pm \equiv \pm 1\ \rm{or} \pm \textit{r} \ ({\rm mod}\ {\textit{b}})$, so $y_2\equiv \pm 1\ \rm{or} \pm \textit{r} \ ({\rm mod}\ {\textit{b}})$. As in the previous case, we obtain a contradiction for $y_2\equiv \pm  r\ ({\rm mod}\ {b})$. From $y_2\equiv \pm  1\ ({\rm mod}\ {b})$, we conclude $y_2=\pm  1$.
\item [\textbullet]  If $W_{2n+1}=q_{2l+1}$, then $rt\equiv  ry_{2}\ ({\rm mod}\ {b})$. Since $\rm{gcd}(\textit{b},\textit{r})=1$, it holds $t\equiv  y_{2}\ ({\rm mod}\ {b})$. Hence, $y_2\equiv \pm 1\ \rm{or} \pm \textit{r} \ ({\rm mod}\ {\textit{b}})$, again. As previous, we have $y_2=\pm  1$.
\end{enumerate}
Therefore, fundamental solutions of the equation (\ref{peta'}) are $(y_2,x_2)=(\pm  1, 1)$.

From the above, we obtain that the following two possibilities are the only types of fundamental solutions that can lead to extensions of the triple $\{2,b,c\}$ to a quadruple $\{2,b,c,d\}$:
\begin{lm}\label{lema-mogucnsti}
For the equation $x=p_l=V_m$, where $p_l$ and $V_m$ are defined by \eqref{osma} and \eqref{deveta}, the following two possibilities exist:
\begin{enumerate}
\item [{a)}] if $l\equiv  m\equiv0\ ({\rm mod}\ {2})$, then $z_0=\pm 1$, $x_0=y_2=1$ and $x_2=1$,
\item [{b)}] if $m\equiv  1\ ({\rm mod}\ {2})$, then $z_0=\pm t$, $x_0=r$, $y_2 =\pm 1$ and $x_2=1$.
\end{enumerate}
\end{lm}

\section{Linear form in three logarithms}

Using techniques from \cite{HPTY}, firstly we transform the equation $p_l=V_m$ into an equality for a linear forms in three logarithms of algebraic numbers. Let \begin{equation}\label{Alpha}\alpha=r+\sqrt{2b}\end{equation} be the solution of Pell equation $t^2-\frac{b}{2}s^2=1$, associated to the Pellian equation (\ref{Pell-D}). Similarly, let \begin{equation}\label{Beta}\beta=s+\sqrt{2c}\end{equation} be the solution of Pell equation $t^2-\frac{c}{2}r^2=1$, associated to the Pellian equation $t^2-\frac{c}{2}r^2=1-\frac{c}{2}$ obtained from (\ref{prva}). Let \begin{equation}\label{gamma}\gamma=\frac{\sqrt{c}(y_2\sqrt{2}+x_2\sqrt{b})}{\sqrt{b}(z_0\sqrt{2}+x_0\sqrt{c})}.\end{equation} We follow the strategy used in \cite{HPTY} and define
\begin{equation}\label{lamda}\Lambda=l\log\alpha-m\log\beta+\log\gamma.\end{equation}As in \cite[Lemma 3]{HPTY}, it can be proven that \begin{equation}\label{nejednakost_za_lamda}0<\Lambda<\beta^{-2m},\end{equation}for $m\geq 1$, which easily leads to $0<\Lambda<\frac{1}{8c}.$ From (\ref{nejednakost_za_lamda}), we have \begin{equation}\label{nejednakost_za_log_lamda}\log\Lambda<-2m\log\beta.\end{equation}


For the possibilies from Lemma \ref{lema-mogucnsti}, we denote
\begin{equation}\label{lambda1}\lambda=\left\{%
\begin{array}{llll}
    0, & \hbox{for\ a),} \\
    1, & \hbox{for\ b)\ with\ $z_0=t$,} \\
    -1, & \hbox{for\ b)\ with\ $z_0=-t$.} \\
\end{array}%
\right.
\end{equation}

In this section we show an upper bound on  $|(l-\lambda)\log\alpha-m\log\beta|$ and we also show that $\Delta = l-\lambda - \nu m \neq 0$. Then we use those results to get a lower bound on index $m$.

\begin{lm} If $x=p_l=V_m$, where $p_l$ and $V_m$ are defined by \eqref{osma} and \eqref{deveta}, then
\begin{equation}\label{aps-lambda}|(l-\lambda)\log\alpha-m\log\beta|=|\Lambda-\lambda\log\alpha-\log\gamma|<\frac{2\sqrt{2}}{\sqrt{b}}.\end{equation}
\end{lm}

\begin{proof} We show this for each possible case in (\ref{lambda1}) separately.
\begin{enumerate}
\item [\textbf{a)}] If $\lambda=0$, then $$\gamma=\frac{\sqrt{c}(\sqrt{2}+\sqrt{b})}{\sqrt{b}(\sqrt{c}\pm \sqrt{2})}=\bigg(1+\frac{\sqrt{2}}{\sqrt{b}}\bigg)\bigg(1+\frac{\sqrt{2}}{\sqrt{c}-\sqrt{2}}\bigg)>1.$$ Hence, $$0<\log\gamma\leq\log  \bigg(1+\frac{\sqrt{2}}{\sqrt{b}}\bigg)+\log\bigg(1+\frac{\sqrt{2}}{\sqrt{c}-\sqrt{2}}\bigg)<\frac{\sqrt{2}}{\sqrt{b}}+\frac{\sqrt{2}}{\sqrt{c}-\sqrt{2}}<\frac{2\sqrt{2}}{\sqrt{b}},$$ which implies (\ref{aps-lambda}).

\item [\textbf{b)}] If $\lambda=1$ and $z_{0}=t$, then $\gamma=\frac{\sqrt{c}(\sqrt{b}\pm\sqrt{2})}{\sqrt{b}(t\sqrt{2}+ r\sqrt{c})}$. In this case $$\alpha^\lambda\gamma-1=\alpha\gamma-1=\frac{\pm\sqrt{2c}(r+\sqrt{2b})-\frac{\sqrt{2b}}{t+\sqrt{bc}}}{\sqrt{b}(t\sqrt{2}+ r\sqrt{c})}.$$ Since $b>4000$, we have $$|\alpha\gamma-1|<\frac{\sqrt{2c}(r+\sqrt{2b})+0.01}{\sqrt{b}(t\sqrt{2}+ r\sqrt{c})}$$ and since $$\sqrt{2c}(r+\sqrt{2b})-\sqrt{2}(t\sqrt{2}+r\sqrt{c})=2\sqrt{bc}-2t<0,$$ it follows that $$|\alpha\gamma-1|<\frac{1.42}{\sqrt{b}}.$$ Furthermore, $$|\log(\alpha\lambda)|=|\log(1+(\alpha\lambda-1))|<\frac{1.42}{\sqrt{b}},$$ so (\ref{aps-lambda}) holds.

\item [\textbf{ c)}]  If $\lambda=-1$ and $z_{0}=-t$, then $\gamma=\frac{\sqrt{c}(\sqrt{b}\pm\sqrt{2})}{\sqrt{b}(-t\sqrt{2}+ r\sqrt{c})}$. In this case $$\alpha^\lambda\gamma-1=\alpha^{-1}\gamma-1=\frac{\sqrt{2}(t\sqrt{bc}-bc)+2\sqrt{b}(r+\sqrt{2b})\pm\sqrt{2c}(t\sqrt{2}+r\sqrt{c})}{\sqrt{b}(r+ \sqrt{2b})(c-2)}.$$ By (\ref{prva}), $t>\sqrt{bc}$, so we have $$|\alpha^{-1}\gamma-1|<\frac{1}{\sqrt{2b}(r+ \sqrt{2b})(c-2)}+\frac{2}{c-2}+\frac{\sqrt{2c}(t\sqrt{2}+r\sqrt{c})}{\sqrt{b}(r+ \sqrt{2b})(c-2)}.$$ Since $b>4000$, from $$\frac{\sqrt{2c}(t\sqrt{2}+r\sqrt{c})}{(r+ \sqrt{2b})(c-2)}-\sqrt{2}=\frac{2\sqrt{2}(r+\sqrt{2b})+\frac{2\sqrt{c}}{t+\sqrt{bc}}}{(r+ \sqrt{2b})(c-2)}<\frac{2\sqrt{2}+0.01}{c-2}<0.01,$$ we get $$|\alpha^{-1}\gamma-1|<\frac{1}{\sqrt{2b}(r+ \sqrt{2b})(c-2)}+\frac{2}{c-2}+\frac{1.42}{\sqrt{b}}<\frac{1.46}{\sqrt{b}}.$$ Finally, $$|\log(\alpha^{-1}\lambda)|=|\log(1+(\alpha^{-1}\lambda-1))|<\frac{1.46}{\sqrt{b}},$$ so (\ref{aps-lambda}) holds.
\end{enumerate}
\end{proof}

As in \cite[Lemma 5]{HPTY}, we prove that the index $l$ is not a multiple of index $m$ increased by $\lambda$, but here the situation is not completely analogue. We have to consider each possible value of $\nu$ separately and apply mathematical induction over $m$. To do that, we will need the following elementary lemma.

\begin{lm} \label{rekurzija}If $(q_m)_{m\geq 0}$ is a second order linear recurrence relation with kernel $(A, B)$, i.~e.~ \begin{eqnarray*} q_{m+2}=Aq_{m+1}+Bq_m,\end{eqnarray*}
then $(q_{2m+1})_{m\geq 0}$ is also a second order linear recurrence relation with kernel $(A^2+2B, -B^2)$,  i.~e.~  \begin{eqnarray*} q_{2m+1} = (A^2+2B)q_{2m-1} - B^2 q_{2m-3}.\end{eqnarray*}
\end{lm}

\begin{proof}
From the initial recurrence relation, we get
\begin{align*}
q_{2m+1} &= A(Aq_{2m-1}+Bq_{2m-2})+Bq_{2m-1} \\
&= A^2q_{2m-1}+B(q_{2m-1}-Bq_{2m-3})+Bq_{2m-1} \\
&=  (A^2+2B)q_{2m-1} - B^2 q_{2m-3}.
\end{align*}
\end{proof}

\begin{lm}\label{lema-Delta} Let $p_l=V_m$, for some positive integers $m$ and $l$, where $p_l$ and $V_m$ are defined by \eqref{osma} and \eqref{deveta}. If the Diophantine quadruple $\{2, b, c, \frac{x^2-1}{2}\}$, where  $c=c_{\nu}^\pm$ for $\nu\in\{1,2,3\}$ and $x=p_l=V_m$,  
is not regular, then
\begin{eqnarray*}\label{lema-Delta1}\Delta = l-\lambda - \nu m \neq 0.\end{eqnarray*}
\end{lm}

\begin{proof}
 Let us define \begin{equation}\label{alpha-nu}\alpha^\nu=(r+\sqrt{2b})^\nu=T_\nu+U_\nu\sqrt{2b},\end{equation} with
\begin{eqnarray}
T_0 &=& 1, \ \  T_1=r, \ \ T_{\nu+2} = 2rT_{\nu+1}-T_\nu,
\label{t_0} \\
U_0 &=& 0, \ \  U_1=1, \ \ U_{\nu+2} = 2rU_{\nu+1}-U_\nu,\label{u_0}\end{eqnarray}for $\nu\geq 0$.  Notice that the sequences $(T_\nu)_{\nu\geq 0}$ and $(U_\nu)_{\nu\geq 0}$ are positive and increasing. From (\ref{t_0}) and (\ref{u_0}), by induction, we get $$T_{m\nu+2\nu}=2T_\nu T_{m\nu+\nu}-T_{m\nu}$$ and $$U_{m\nu+2\nu}=2T_\nu U_{m\nu+\nu}-U_{m\nu}.$$ Hence, from (\ref{deveta}),
$$p_{m\nu+2\nu}=2T_\nu p_{m\nu+\nu} - p_{m\nu}.$$ Also, it is easy to show that $p_\nu = x_2T_\nu+2y_2 U_\nu$. By 
(\ref{Pell-general}), (\ref{rekurzija1}), (\ref{t_0}) and (\ref{u_0}), we have \begin{equation}\label{s}s=s_\nu^\pm=T_{\nu}\pm 2U_\nu,\end{equation}for $\nu\in\{1,2,3\}$. We distinguish three cases, depending on the value of $\lambda$ in (\ref{lambda1}).
\begin{enumerate}
\item [\textbf{a)}] If $\lambda=0$, then $z_0=\pm 1, x_0=1, y_2=x_2=1$ so $p_\nu =  T_\nu + 2U_\nu$. In (\ref{osma}) and (\ref{deveta}), we have
\begin{eqnarray*} 
V_0 &=& 1, \ \  V_1=s\pm 2, \ \ V_{m+2} = 2sV_{m+1}-V_m, \\
p_0 &=& 1, \ \  p_1=r+2, \ \ p_{l+2} = 2rp_{l+1}-p_l,\end{eqnarray*}for $m,l\geq 1$. If we prove that $V_m \neq p_{m\nu}$, then $l\neq m\nu$ so, in this case, $\Delta\neq 0$. More precisely, by induction over $m$ we will prove that if $\{a, b, c, d\}$ is irregular, then $V_m \neq p_{m\nu}$. We distinguish two cases in (\ref{s}). Notice that in this part of the proof we consider the case from Lemma \ref{lema-mogucnsti} a), so $x=V_m$ is not even possible for odd values of $m$. 

\begin{enumerate}
\item [\textbf{a1)}] Let $s=s_\nu^+=T_\nu+2U_\nu$.

Let $V_1=s-2$. For $m=1$, we have $V_1 < p_\nu$. For $m=2$ and $\nu = 1$, we have $x=V_2=p_2=2s(s-2)-1$. In that case, we get a regular quadruple $\{a, b, c, d\}$, where $c=2k^2+6k+4$ and  $d=d_{+}=4 (3 + 20 k + 42 k^2 + 32 k^3 + 8 k^4)$.  For $\nu > 1$, $$V_2=2s(s-2)-1 > 2T_\nu p_\nu-p_0 =p_{2\nu}$$ is equivalent to
\begin{eqnarray*}
  (T_\nu+2U_\nu)(T_\nu+2U_\nu - 2) > T_\nu (T_\nu+2U_\nu) 
&\Leftrightarrow & T_\nu+2U_\nu - 2 > T_\nu\\
&\Leftrightarrow & U_\nu > 1.
\end{eqnarray*}For $m=3$ and for $\nu\geq 1$, we have $$V_3  =2sV_2 - (s-2) > 2T_\nu p_{2\nu}-p_\nu = p_{3\nu},$$ because this is equivalent with 
\begin{eqnarray*}
& &2s(2s(s-2)-1) - T_\nu-2U_\nu+2 > 2T_\nu(2T_\nu s-1)-T_\nu-2U_\nu \\
& &\Leftrightarrow  2s(2s^2-4s) - 2(T_\nu+2U_\nu) +2> 4T_\nu^2 s-2T_\nu\\
& & \Leftrightarrow  2s(s^2-2s-T_\nu^2)-2U_\nu +1 > 0,
\end{eqnarray*}
which is true since $s^2-2s-T_\nu^2 > 1$.


If $V_1=s+2$, then $V_1 > p_\nu$ and $V_2 > p_{2\nu}$ for $\nu\geq 1$.\\

In both subcases $V_1=s-2$  and $V_1=s+2$ of the case \textbf{a1)} the step of induction is the same. We assume $p_{m\nu+\nu} < V_{m+1}$ which implies
\begin{align*}
p_{m\nu + 2\nu} = 2T_\nu p_{m\nu+\nu} - p_{m\nu}
&< 2 T_\nu V_{m+1} + (4U_\nu V_{m+1}-V_m) \\
& =2(T_\nu+2U_\nu)V_{m+1}-V_m \\
&= 2sV_{m+1}-V_m \\
&= V_{m+2}.
\end{align*}
Hence, we have proven $V_m\neq p_{m\nu}$ in the case \textbf{a1)}. 

\item [\textbf{a2)}] Let $s=s_\nu^-=T_\nu-2U_\nu$.

For $\nu=1$, we have $s=s_1^-=r-2$, which is not possible. For $\nu \in  \{2,3 \}$, we have the following situatuation. For $m=1$, we have $$p_\nu = T_\nu+2U_\nu > T_\nu-2U_\nu \pm 2=V_1,$$  which is equivalent to $U_\nu > \pm \frac{1}{2}$.

If $m=2$, there holds
\begin{eqnarray*}
V_2 < p_{2\nu} 
&\Leftrightarrow&  2s^2\pm 4s - 1 < 2T_\nu^2+4T_\nu U_\nu-1 \\
&\Leftrightarrow&  2(T_\nu^2-4T_\nu U_\nu+4U_\nu^2)\pm 4(T_\nu-2U_\nu)-1 < 2T_\nu^2+4T_\nu U_\nu-1\\
&\Leftrightarrow&  2U_\nu^2\pm (T_\nu-2U_\nu) < 3T_\nu U_\nu,
\end{eqnarray*}
which is true because $T_\nu > 2U_\nu$.\\

Before doing step of induction, let us notice that $2T_\nu-1 > 2s$ and that $(p_{m\nu})_{m\geq 0}$ is increasing sequence. Assume now that $p_{m\nu+\nu} > V_{m+1}$. Then, \begin{eqnarray*}p_{m\nu+2\nu} &=& 2T_\nu p_{m\nu+\nu} - p_{m\nu} \\
&=& (2T_\nu-1)p_{m\nu+\nu} + p_{m\nu+\nu}-p_{m\nu} \\
&>& 2sV_{m+1}\\
&>& 2sV_{m+1}-V_m\\
&=&V_{m+2}.
\end{eqnarray*}Therefore, $V_m\neq p_{m\nu}$ in the case \textbf{a2)}. 
\end{enumerate}

\item [\textbf{b)}] If $\lambda=1$, then $z_0=t, x_0=r, x_2=1, y_2=\pm 1$ so $p_\nu =  T_\nu \pm 2U_\nu$. From (\ref{osma}) and (\ref{deveta}), we have:
\begin{eqnarray*}
V_0 &=& r, \ \  V_1=rs+2t, \ \ V_{m+2} = 2sV_{m+1}-V_m,\\
p_0 &=& 1, \ \  p_1=r\pm 2, \ \ p_{l+2} = 2rp_{l+1}-p_l,\end{eqnarray*}for $m,l\geq 1$. Notice that $(V_m)_{m\geq 0}$ and $(p_l)_{l\geq 0}$ are increasing sequences.

If we prove that $V_m \neq p_{m\nu+1}$, then $l\neq m\nu+1$ so $\Delta\neq 0$. By induction over $m$ we will prove that if $\{a, b, c, d\}$ is irregular, then $V_m \neq p_{m\nu+1}$. Notice that here we consider the case from Lemma \ref{lema-mogucnsti} b), so $x=V_m$ is not possible at all for even values of $m$. 

For $m=1$, if $x=V_1=p_{\nu+1}$ then, by (\ref{d-st-reg}),  $d=d_{+}$. For $m\geq 2$, the idea is to express $s$ and $t$ as polynomials in variable $r$. Hence, we also have $V_m$ and $p_l$ expressed as polynomials in variable $r$ and we can compare them. For $\nu=1$, $s=s_1^{+}=r+2$ and $t= t_1^{+}=(r^2+2r-1)/2$. We obtain $p_3<V_2 $ and 
$p_4<V_3$. For $\nu=2$, if $s=s_2^{-}=2r^2-4r-1$ and $t=t_2^{-}=r^3-2r^2-r+1$, it holds $p_5>V_2$ and $p_7>V_3$. 
Similarly, if $s=s_2^{+}=2r^2+4r-1$, then $p_5<V_2$ and $p_7<V_3$. For $\nu=3$, if $s=s_3^{-}=4r^3-8r^2-3r+2$ then $p_7>V_2$ and $p_{10}>V_3$, while if $s=s_3^{+}=4r^3+8r^2-3r-2$ then $p_7<V_2$ and $p_{10}<V_3$. Therefore, for $m=2$ and $m=3$, we have $V_m\neq p_{\nu m+1}$.\\ 

The step of induction can be obtained for each possible value of $\nu$ separately, as follows. For $\nu=1$, if $p_m<V_{m-1}$ then \begin{eqnarray*}p_{m+1}=2rp_{m}-p_{m-1}&<&2rV_{m-1}+2V_{m-1}-V_{m-2}\\
&=&2sV_{m-1}-V_{m-2}\\
&=&V_{m}.\end{eqnarray*}

For $\nu=2$, we will use Lemma \ref{rekurzija}.

Let $s=s_2^{-}$. If $p_{2m-1}>V_{m-1}$ then, by Lemma \ref{rekurzija}, \begin{eqnarray*}p_{2m+1}=(4r^2-2)p_{2m-1}-p_{2m-3}&>&(4r^2-3)p_{2m-1}\\
&>&(4r^2-3)V_{m-1}\\
&>&\frac{4r^2-3}{2s}V_m\\
&=&\frac{4r^2-3}{4r^2-8r-2}V_m>V_m.\end{eqnarray*}

Let $s=s_2^{+}$. If $p_{2m-1}<V_{m-1}$ then, by Lemma \ref{rekurzija}, \begin{eqnarray*}p_{2m+1}=(4r^2-2)p_{2m-1}-p_{2m-3}&<&(4r^2-2)p_{2m-1}\\
&<&(4r^2-2)V_{m-1}\\
&<&\frac{4r^2-2}{2s-1}V_m\\
&=&\frac{4r^2-2}{4r^2+8r-3}V_m<V_m.\end{eqnarray*}

For $\nu=3$, we firstly get
\begin{eqnarray*}p_{3m+1}&=&2rp_{3m}-p_{3m-1}\\
&=&2r(2rp_{3m-1}-p_{3m-2})-p_{3m-1}\\
&=&(4r^2-1)p_{3m-1}-2rp_{3m-2}\\
&=&(4r^2-1)(2rp_{3m-2}-p_{3m-3})-2rp_{3m-2}\\
&=&(8r^3-4r)p_{3m-2}-(4r^2-1)p_{3m-3}.
\end{eqnarray*}Hence, \begin{equation}\label{donja_i_gornja}(8r^3-4r^2-4r+1)p_{3m-2}<p_{3m+1}<(8r^3-4r)p_{3m-2}.\end{equation}

Let $s=s_3^{-}$. If $p_{3m-2}>V_{m-1}$ then, by (\ref{donja_i_gornja}), \begin{eqnarray*}p_{3m+1}>(8r^3-4r^2-4r+1)p_{3m-2}
&>&(8r^3-4r^2-4r+1)V_{m-1}\\
&>&\frac{8r^3-4r^2-4r+1}{2s}V_m\\
&=&\frac{8r^3-4r^2-4r+1}{8r^3-16r^2-6r+4}V_m>V_m.\end{eqnarray*}

Let $s=s_3^{+}$. If $p_{3m-2}<V_{m-1}$ then, by (\ref{donja_i_gornja}), \begin{eqnarray*}p_{3m+1}<(8r^3-4r)p_{3m-2}
&<&(8r^3-4r)V_{m-1}\\
&<&\frac{8r^3-4r}{2s-1}V_m\\
&=&\frac{8r^3-4r}{8r^3+16r^2-6r-3}V_m<V_m.\end{eqnarray*}

\item [\textbf{c)}] We omit the proof as it is very similar to the proof of case \textbf{b)}.

\end{enumerate}
\end{proof}
\begin{remark}If $l=0$ or $m=0$ in equation $x=p_l=V_m$, then $d=0$ by (\ref{treca}), \eqref{osma}, \eqref{deveta} and Lemma \ref{lema-mogucnsti}. If $m=1$, then we can obtain a regular quadruple  $\{2, b, c, \frac{x^2-1}{2}\}$, as it is explained in parts \textbf{b)} and \textbf{c)} of the Lemma  \ref{lema-Delta}. Since in the part \textbf{a)} the case $x=V_1$ is not possible, for $m=1$ the quadruple $\{2, b, c, \frac{x^2-1}{2}\}$ is regular, if it exists. If $m=2$, then we can also obtain a regular quadruple  $\{2, b, c, \frac{x^2-1}{2}\}$, as it is explained in the part \textbf{a1)} of the Lemmma  \ref{lema-Delta}.  By Fujita \cite[Lemma 8]{fujita}, for $m=2$ we can't have an irregular Diophantine quadruple. 
\end{remark}

Now we want to find the lower bound for $m$ in a solution $(l,m)$ of the equation $p_l=V_m$. 
\begin{lm}\label{lema4}If the equation $p_l=V_m$, where $p_l$ and $V_m$ are defined by \eqref{osma} and \eqref{deveta},  has a solution $(l,m)$ with $m\geq 1$ then, for $b>4000$, we have \begin{equation}\label{m}m>0.69|\Delta|\sqrt{b}\log\alpha.\end{equation}
\end{lm}

\begin{proof}
From (\ref{aps-lambda}), we obtain \begin{eqnarray*}\bigg|\frac{l-\lambda}{m}-\frac{\log\beta}{\log\alpha}\bigg|<\frac{2\sqrt{2}}{m\sqrt{b}\log\alpha}.\end{eqnarray*} Using that and (\ref{lema-Delta1}), we further obtain \begin{equation}\label{Delta_m}\frac{|\Delta|}{m}<\bigg|\frac{\log\beta}{\log\alpha}-\nu\bigg|+\frac{2\sqrt{2}}{m\sqrt{b}\log\alpha}.\end{equation} Also, there holds \begin{equation}\label{logaritmi}\bigg|\frac{\log\beta}{\log\alpha}-\nu\bigg|=\bigg|\frac{\log(\frac{\beta}{\alpha^\nu})}{\log\alpha}\bigg|=\bigg|\frac{\log(1+\frac{\beta-\alpha^\nu}{\alpha^\nu})}{\log\alpha}\bigg|.\end{equation}Using (\ref{Beta}), (\ref{alpha-nu}) and (\ref{s}), we have \begin{eqnarray*}\bigg|\frac{\beta-\alpha^\nu}{\alpha^\nu}\bigg|&=&\bigg|\frac{s+\sqrt{2c}-(r+\sqrt{2b})^\nu}{(r+\sqrt{2b})^\nu}\bigg|=\Bigg|\frac{2s-\frac{1}{s+\sqrt{2c}}-\Big(2T_\nu-\frac{1}{T_\nu+U_\nu\sqrt{2b}}\Big)}{T_\nu+U_\nu\sqrt{2b}}\Bigg|\\
&=&\bigg|\frac{\pm 4U_\nu-\frac{1}{s+\sqrt{2c}}+\frac{1}{T_\nu+U_\nu\sqrt{2b}}}{T_\nu+U_\nu\sqrt{2b}}\bigg|<\frac{4U_\nu+0.01}{2U_\nu\sqrt{2b}}<\frac{1.42}{\sqrt{b}}.\end{eqnarray*}Hence, \begin{equation}\label{ograda-log}\bigg|\log\Big(1+\frac{\beta-\alpha^\nu}{\alpha^\nu}\Big)\bigg|<1.01\bigg|\frac{\beta-\alpha^\nu}{\alpha^\nu}\bigg|<\frac{1.44}{\sqrt{b}}.\end{equation}From (\ref{Delta_m}), (\ref{logaritmi}) and (\ref{ograda-log}),
\begin{eqnarray*}\frac{|\Delta|}{m}<\frac{1.44}{\sqrt{b}\log\alpha}+\frac{2\sqrt{2}}{m\sqrt{b}\log\alpha}=\frac{1.44+\frac{2\sqrt{2}}{m}}{\sqrt{b}\log\alpha}\end{eqnarray*}and then
\begin{eqnarray*}1.44m+2\sqrt{2}>|\Delta|\sqrt{b}\log\alpha.\end{eqnarray*}Therefore, (\ref{m}) holds.
\end{proof}


\section{Linear forms in two logarithms and the proof of Theorem \ref{Tm1}}

As in \cite{HPTY}, we apply Laurent's results (see \cite{L}) on linear forms in two logarithms. We obtain an upper bound on $m$, which will contradict lower bound from Lemma \ref{lema4} unless $k < 10^6$ roughly. Then we finish the proof of our main result using the well known Baker-Davenport reduction method.


We can rewrite (\ref{lamda}) into
\begin{equation}\label{lamda1}\Lambda=\log(\alpha^{\Delta+\lambda}\gamma)-m\log\Big(\frac{\beta}{\alpha^\nu}\Big)=m\log\Big(\frac{\alpha^\nu}{\beta}\Big)-\log(\alpha^{-\Delta-\lambda}\gamma^{-1}),\end{equation}where $\alpha$, $\beta$, $\gamma$ and $\alpha^\nu$ are given with (\ref{Alpha}), (\ref{Beta}), (\ref{gamma}) and (\ref{alpha-nu}), respectively. Let $\alpha_1:=\frac{\alpha^\nu}{\beta}$ and $\alpha_2:=\alpha^{\Delta+\lambda}\gamma$. Then, $\alpha_1$ is a zero of the polynomial
\begin{eqnarray*}\bigg(X-\frac{T_\nu+U_\nu\sqrt{2b}}{s+\sqrt{2c}}\bigg)\bigg(X-\frac{T_\nu-U_\nu\sqrt{2b}}{s+\sqrt{2c}}\bigg)\bigg(X-\frac{T_\nu+U_\nu\sqrt{2b}}{s-\sqrt{2c}}\bigg)\bigg(X-\frac{T_\nu-U_\nu\sqrt{2b}}{s-\sqrt{2c}}\bigg)\end{eqnarray*}
$$=X^4-4sT_\nu X^3+(4T_\nu^2+8c+1)X^2-4sT_\nu X+1,$$ which is its minimal polynomial over $\mathbb{Z}$, or minimal polynomial divides it.

For any non-zero algebraic number $\alpha$ of degree $d$ over $\mathbb{Q}$, with minimal polynomial $a_0\prod_{j=1}^{d}(X-\alpha^{(j)})$, the absolute logarithmic height of $\alpha$ is defined with $$h(\alpha)=\frac{1}{d}\bigg(\log |a_0|+\sum_{j=1}^d \log\max \{1,|\alpha^{(j)}|\} \bigg),$$where $\alpha^{(j)}$ are the conjugates of $\alpha$ in $\mathbb{C}$. Here we have the linear form (\ref{lamda1}) in two algebraic numbers $\alpha_1$ and $\alpha_2$ over $\mathbb{Q}$. Since at most two conjugates of $\alpha_1$ are greater than $1$, depending on wether $\alpha^\nu>\beta$ or $\alpha^\nu<\beta$, as in \cite{HPTY}, we have $$h(\alpha_1)\leq\frac{\nu}{2}\log\alpha\ \ \rm{or}\ \  \textit{h}(\alpha_1)\leq\frac{1}{2}\log\beta.$$

It holds $$h(\alpha^{\Delta+\lambda})=\frac{1}{2}|\Delta+\lambda|\log\alpha$$ and
\begin{eqnarray*}h(\gamma)=h\bigg(\frac{\sqrt{c}(y_2\sqrt{2}+x_2\sqrt{b})}{\sqrt{b}(z_0\sqrt{2}+x_0\sqrt{c})}\bigg)&\leq& h\bigg(\frac{y_2\sqrt{2}+x_2\sqrt{b}}{\sqrt{b}}\bigg)+h\bigg(\frac{z_0\sqrt{2}+x_0\sqrt{c}}{\sqrt{c}}\bigg)\\
&\leq& \frac{1}{2}\log (b+\sqrt{2b})+\frac{1}{2}\log (rc+t\sqrt{2c})\\
&<& \frac{1}{2}\log (4rbc)<\frac{3}{2}\log \alpha+\log \beta.
\end{eqnarray*}Therefore, \begin{eqnarray*}h(\alpha_2)=h(\alpha^{\Delta+\lambda}\gamma)\leq \frac{1}{2}(|\Delta+\lambda|+3)\log\alpha+\log\beta.\end{eqnarray*}By (\ref{ograda-log}), $\frac{|\alpha^\nu-\beta|}{\alpha^\nu}<\frac{1.426}{\sqrt{b}}$. Assuming $k\geq1000$, we have $$|\log\alpha^\nu-\beta|<0.01.$$ We use the notation as in \cite[Lemma 8]{HPTY} and get \begin{eqnarray*}h_1&=&\frac{\nu}{2}\log\alpha+0.01>h(\alpha_1),\\
h_2&=& \frac{1}{2}(|\Delta+\lambda|+3+2\nu)\log\alpha+0.01>h(\alpha_2).\end{eqnarray*}Further, since by (\ref{lamda1}) $b_1=m$, $b_2=1$ and $D=4$, we have $$\frac{|b_2|}{Dh_1}=\frac{1}{2\nu\log\alpha+0.04}<0.07.$$ Let us define \begin{equation}\label{b_crtica}b'=\frac{m}{2(|\Delta+\lambda|+3+2\nu)\log\alpha+0.04}+0.07.\end{equation}If $\log b'+0.38\leq\frac{30}{D}=7.5$, then $$b'\leq 1236.$$ Else, by \cite[Lemma 8]{HPTY},
\begin{eqnarray*}\log |\Lambda|\geq -17.9\cdot 4^4 (\log b'+0.38)^2 h_1h_2.\end{eqnarray*}

Also, by (\ref{nejednakost_za_log_lamda}), $$m\log\beta<17.9\cdot 128 (\log b'+0.38)^2 h_1h_2.$$ Since $\log\beta>\log\alpha^\nu-0.01>2h_1-0.03$, we have $$m<1.01\cdot 17.9\cdot 64 (\log b'+0.38)^2 h_2$$ and then $$b'-0.07=\frac{m}{4h_2}<289.264(\log b'+0.38)^2,$$ which yields $b'\leq 33789.$

Finally, from (\ref{b_crtica}), we obtain the following statement. 
\begin{lm}\label{lema5}If for the triple $\{2,b,c_\nu^\pm\}$, with $1\leq \nu \leq 3$, the equation $p_l=V_m$, where $p_l$ and $V_m$ are defined by \eqref{osma} and \eqref{deveta}, has a solution $(l,m)$ with $m\geq 1$ then, for $k\geq 1000$, we have \begin{equation}\label{m_gornja}m<67578.15(|\Delta+\lambda|+3+2\nu)\log\alpha+1351.57.\end{equation}\end{lm}

\bigskip

If the Diophantine quadruple $\{2, b, c, \frac{x^2-1}{2}\}$, where  $c=c_{\nu}^\pm$ for $\nu\in\{1,2,3\}$ and $x=p_l=V_m$,  
is not regular, then by Lemma \ref{lema-Delta}, $\Delta\neq 0$. Therefore, we assume that $|\Delta|\geq 1$. If $k\geq 1000$, by Lemma \ref{lema4} and Lemma \ref{lema5}, we get $$0.69|\Delta|\sqrt{b}\log\alpha<67578.15(|\Delta+\lambda|+3+2\nu)\log\alpha+1351.57.$$ From that, we have
\begin{eqnarray*}\sqrt{b}&<&\frac{67578.15(|\Delta+\lambda|+3+2\nu)}{0.69|\Delta|}+\frac{1351.57}{0.69|\Delta|\log\alpha}\\
&<&97939.36(2\nu+5)+257.71\\
&<&1077590.56.\end{eqnarray*}By inserting (\ref{k}) into the last inequality, we obtain $$k\leq 761970.$$

We now finish the proof of Theorem \ref{Tm1} using the Baker-Daveport reduction method, which is standard method in solving such problems for years. We use its version from \cite{dujepet}. We get the first bound $m<1.33\cdot 10^{18}$ using the same method as it was 
used in \cite[ Section 8]{duje1}. Only this time we have the exact values for 
fundamental solutions. In at most two steps of reduction in all cases we get $m\leq3$. By Fujita \cite[Lemma 8]{fujita}, if the equation (\ref{jednadzbica}) has a solution which leads to an irregular Diophantnine quadruple, then $m,n\geq 3$ and $(m,n) \neq 3$. Hence, if $m=3$ then $n\geq 4$. By Dujella \cite[Lemma 3]{duje1}, for $m=3$ it holds $n\leq 4$. Therefore, $n=4$, which is not possible since $m$ and $n$ have the same parity. That finishes the proof of Theorem \ref{Tm1} since for $m\leq2$, by Remark 1,  we get only the extensions to a regular quadruple (or $d=0$).  

\bigskip

\textbf{Acknowledgement:} The authors were supported by the Croatian Science Foundation under the project no. IP-2018-01-1313.

\bigskip
Faculty of Civil Engineering, University of Zagreb,\\ Fra Andrije Ka\v{c}i\'{c}a-Mio\v{s}i\'{c}a 26, 10000 Zagreb, Croatia \\
Email: nadzaga@grad.hr \\
\\
Faculty of Civil Engineering, University of Zagreb,\\ Fra Andrije Ka\v{c}i\'{c}a-Mio\v{s}i\'{c}a 26, 10000 Zagreb, Croatia \\
Email: filipin@grad.hr \\
\\
Department of Mathematics, University of
Rijeka,\\ Radmile Matej\v{c}i\'{c} 2, 51000 Rijeka, Croatia \\
Email: ajurasic@math.uniri.hr

\end{document}